\documentclass[12pt]{article}
\usepackage{amsfonts}

\usepackage{latexsym}
\usepackage{amssymb}
\usepackage{epsfig}
\usepackage{amsmath}
\usepackage{amsthm}
\usepackage[mathscr]{eucal}
\usepackage{subfigure}
\usepackage{graphicx}
\usepackage{hyperref}

\newtheorem{Lemma}{Lemma}

\newtheorem{Theorem}[Lemma]{Theorem}

\newtheorem{Definition}{Definition}

\renewcommand{\qed}{\hfill{\ \ \rule{2mm}{2mm}} \vspace{0.2in}}

\newcommand{\ind}{1\hspace{-2.3mm}{1}}

\begin{document}
\title{Weighted Domination and Colouring in Random Graphs}
\author{ \textbf{Ghurumuruhan Ganesan}
\thanks{E-Mail: \texttt{gganesan82@gmail.com} } \\
\ \\
IISER, Bhopal}
\date{}
\maketitle

\begin{abstract}
In the first part of this paper, we consider weighted domination in the case where the vertices of the complete graph on~\(n\) vertices are equipped with independent and identically distributed (i.i.d.) weights. We use the probabilistic iteration to determine sufficient conditions for maximizing the weighted domination probability. In the second part, we study  a weighted generalization of the chromatic number and estimate the minimum number of colours needed to satisfy the constraints when the weights themselves are random. We show that the ``extra" cost incurred for weighted colouring is small if the weights have sufficiently large moments. We also consider inhomogenous random graphs where the edge probabilities are not necessarily all same and obtain bounds for the chromatic number in terms of its averaged edge probabilities. 


\vspace{0.1in} \noindent \textbf{Key words:} Weighted Domination; Maximum Domination Probability; Weighted Colouring; Inhomogenous Random graphs,

\vspace{0.1in} \noindent \textbf{AMS 2020 Subject Classification:} Primary: 60C05.
\end{abstract}

\bigskip

\setcounter{equation}{0}
\renewcommand\theequation{\thesection.\arabic{equation}}
\section{Introduction}
In this section, we describe briefly the two problems studied in this paper.

\subsection*{Weighted Domination}
Records in a sequence of independent and identically distributed (i.i.d.) random variables have been extensively studied in the context of statistical estimation. The probability that the~\(k^{th}\) random variable is a record (as compared to the previous~\(k-1\) values) grows inversely in the index~\(k\) and the corresponding record events are mutually independent~\cite{nev}. Various aspects of record properties with differing assumptions have been studied since and as examples,~\cite{bair} study large deviation properties for weak records and~\cite{jan} studies characterizations of geometric distributions based on~\(k^{th}\) record values. Also,~\cite{chara} study records in a sequence of maximally dependent random variables and obtain expressions for corresponding record moments and distributions.

In this paper, we study records from a graph theoretic view point and use probabilistic iteration to estimate the probability that a given set of vertices form a weighted dominating set. 

\subsection*{Weighted Colouring}
The chromatic number~\(\chi(G)\) of a graph~\(G\) is the smallest number of colours needed for a proper colouring of~\(G\) and many upper and lower bounds for~\(\chi(G)\) exist in terms of various graph parameters like maximum vertex degree, independence number, clique number etc. (see Chapter~\(5,\)~\cite{west}). For homogenous random graphs where each edge is present with the same probability~\(p_n,\) concentration of the chromatic number around its mean has been well-studied~\cite{sham,boll,sct}. Many bounds for the expected chromatic number have also been derived using clique sizes in the complement graph and a combination of second moment and martingale methods~\cite{bol2,alon2}. The paper~\cite{pan} uses counting arguments to obtain sharp lower bounds on the chromatic number. In a related direction,~\cite{alon3} obtain two point concentration of the chromatic number for edge probability~\(p_n = n^{-1/2-\delta}\) with~\(\delta > 0.\) The analysis proceeds through~\(k-\)choosable graphs. Extending this to~\(p_n = \frac{d}{n},\)~\cite{ach} use analytical techniques to obtain the two possible values of the chromatic number for~\(d > 0\) a constant.

Recently, in~\cite{gan} we have studied the chromatic number of homogenous random graphs whose edge probabilities do not necessarily form a convergent sequence. In the first part of this paper, we obtain estimates for the chromatic number of an inhomogenous graph in terms of averaged edge probabilities. We use an auxiliary bound for chromatic number of a deterministic graph, in terms of its maximum averaged degree, that is of independent interest. 

In the second part of this paper, we study weighted colouring in random graphs. Generalizing radio labelling of graphs~\cite{chart,korz} which we call as weighted colouring number, we use the probabilistic method to estimate the weighted colouring number when the edge weights themselves are random. Such situations arise frequently in applications and we show that if the weights have sufficiently large moments, then the ``cost" incurred due to the unboundedness of weights is small.

The paper is organized as follows: In the following section, we state and prove our first main result regarding maximizing the domination probability in weighted random graphs. Next, in Section~\ref{sec_col} we state and prove our two main results regarding weighted and constrained colouring in random graphs and finally, in Section~\ref{sec_comb_lem}, we collect the auxiliary combinatorial lemmas used in the proof of the main theorems.

\setcounter{equation}{0}
\renewcommand\theequation{\thesection.\arabic{equation}}
\section{Weighted Domination}\label{sec_misc}
In this section, we study domination in weighted random graphs obtained as follows. Let~\(K_n\) be the complete graph on~\(n\) vertices and equip vertex~\(v\) with a random weight~\(S_v.\) The random variables~\(\{S_v\}_{1 \leq v \leq n}\) are i.i.d.\ continuous with a common cumulative distribution function (cdf)~\(F.\) For each vertex~\(v,\) we now associate a \emph{deterministic} domination neighbourhood~\({\cal D}(v) \subset \{1,2,\ldots,n\} \setminus \{v\}\) and say that~\(v\) is a \emph{weighted dominating vertex} if~\(S_v > \max_{u \in {\cal D}(v)} S_u.\) Similarly, we say that a finite set of vertices~\({\cal V} := \{v_1,\ldots,v_k\}\) is a \emph{weighted dominating set} if each vertex in~\({\cal V}\) is a weighted dominating vertex.

Let~\(E_{dom} = E_{dom}({\cal V})\) be the event that~\({\cal V}\) is a weighted dominating set and also let~\(R_{tot} := \sum_{1 \leq l \leq k} \ind(A_{v_l})\) be the total number of weighted dominating vertices in~\({\cal V}.\) For any vertex~\(v,\) we have by symmetry that~\(\mathbb{P}(A_v) = \frac{1}{d(v)},\) where~\(d(v)=1+\#{\cal D}(v)\) and~\(\#{\cal D}(v)\) is the size of the domination neighbourhood of~\(v.\) Therefore if the vertices in~\({\cal V}\) have disjoint domination neighbourhoods, then the events~\(A_{v}, v \in {\cal V}\) are independent and so~\[\mathbb{P}(E_{dom}) = \prod_{l=1}^{k} \frac{1}{d(v_l)} =: p_{dom}.\]

In general, if the domination neighbourhoods are not disjoint, the events~\(A_{v}, v \in {\cal V}\) are correlated and so we expect that~\(E_{dom}\) occurs with lesser probability, as described in the following result. Recalling that~\({\cal V} = \{v_1,\ldots,v_k\}\) with~\(v_1<v_2 < \ldots < v_k,\) we define~\({\cal C}(v_l) := \bigcup_{j=1}^{l} {\cal D}(v_j)\) and set~\(c(v_l) := 1+ \#{\cal C}(v_l)\) for~\(1 \leq l \leq k.\)
\begin{Theorem}\label{prop1} If~\(\int_{0}^{x} F^{k-1}(y) dF(y) = \frac{F^{k}(x)}{k}\) for each integer~\(k \geq 1\) and~\(x > 0\)
and
\begin{equation}\label{inclu_cond}
v_k \notin {\cal C}(v_k) \text{ and } v_l \in {\cal D}(v_{l+1}) \setminus {\cal C}(v_l) \text{ for each } 1 \leq l \leq k-1,
\end{equation}
then
\begin{equation}\label{mut_ind}
\mathbb{P}\left(E_{dom}\right) = \prod_{1 \leq l \leq k} \frac{1}{c(v_l)} \leq p_{dom} \text{ and }var(R_{tot}) \leq \mathbb{E}R_{tot} = \sum_{l=1}^{k} \frac{1}{d(v_l)}.
\end{equation}
Moreover, the events~\(\{A_{v_l}\}_{1 \leq l \leq k}\) are mutually independent if and only if\\\({\cal D}(v_l) \subset {\cal D}(v_{l+1})\) strictly, for each~\(1 \leq l \leq k-1.\)
\end{Theorem}
In words, the above result says that if the domination neighbourhoods satisfy the ``weak-nested" property~(\ref{inclu_cond}), then~\(p_{dom}\) is the maximum possible weighted domination probability and moreover, this value is achieved only if the neighbourhoods satisfy a strong-nested property.

From a statistical view point, we could also interpret Theorem~\ref{prop1} as a ``graph-theoretic" version of records. For context, we recall that ``time-based" record events in a sequence of i.i.d.\ random variables are known to be mutually independent when the comparison set consists of the entire past~\cite{nev}. From Theorem~\ref{prop1}, we see that graph-theoretic record events are \emph{negatively correlated} for weakly nested neighbourhoods and are mutually independent for strongly nested neighbourhoods.











\emph{Proof of Theorem~\ref{prop1}}: Let~\(B_k := \bigcap_{j=1}^{k} A_{v_j} \) and for~\(x > 0\)   define
\begin{equation}\label{akx}
A_{v_k}(x) := A_{v_k} \bigcap \left\{S_{v_k} < x\right\} = \left\{ \max_{j \in {\cal D}(v_k)} S_j < S_{v_k} < x\right\}
\end{equation}
 and~\(B_k(x) := B_{k-1} \bigcap A_{v_k}(x),\) with the notation that the maximum of the empty set is zero and~\(B_0 = \Omega.\)  We first show by induction that for every~\(x > 0,\)
\begin{equation}\label{homo_b}
\mathbb{P}(B_k(x)) = \frac{F^{c(v_k)}(x)}{c(v_1)\cdots c(v_k)} = \mathbb{P}(B_k) \cdot F^{c(v_k)}(x).
\end{equation}
For the basis step, we see that the random variable~\(X\) with cdf~\(F\) is continuous and so~\(F(x) = \mathbb{P}(X \leq x) = \mathbb{P}(X <x).\) Conditioning on~\(S_{v_1} = y,\) we therefore get that~\(\mathbb{P}\left(B_1(x)\right) = \int_{0}^{x} \mathbb{P}\left(\max_{j \in {\cal D}(v_1)} S_j < y\right)  dF(y)\) equals
\begin{equation}
\int_{0}^{x} \mathbb{P}\left(\bigcap_{j \in {\cal D}(v_1)} \{S_j < y\}\right) dF(y)  = \int_{0}^{x} F^{d(v_1)-1}(y) dF(y) = \frac{F^{d(v_1)}(x)}{d(v_1)}.\label{pr_ak_ev}
\end{equation}
Since~\(d(v_1) = c(v_1)\) this proves the basis step.

To prove the induction step, we now assume that the relation~(\ref{homo_b}) is true for the event~\(B_{k-1}(x)\) and consider the event~\(B_{k}(x).\) If~\(B_{k-1} = \bigcap_{j=1}^{k-1} A_{v_j}\) occurs, then using the weak nested property in the statement of the theorem, we already have~\(S_{v_{k-1}} > \max_{l \in {\cal C}(v_{k-1})}S_{v_l}\) and so~\(S_{v_k}\) is a record if and only if~\(S_{v_k} > \max_{j \in {\cal E}_k} S_j\) and~\(S_{v_k} > S_{v_{k-1}},\) where~\({\cal E}_k := {\cal D}(v_k) \setminus \left({\cal C}(v_{k-1}) \cup v_{k-1}\right).\) Letting~\(\ind(.)\) denote the indicator function, we then have that
\begin{eqnarray}
\ind\left(B_k(x)\right) &=& \ind\left(A_{v_k}(x) \bigcap B_{k-1}\right) \nonumber\\
&=& \ind\left(x > S_{v_k} > \max_{j \in {\cal E}_k}  S_j\right) \ind(x > S_{v_k} > S_{v_{k-1}}) \ind(B_{k-1}) \label{ai2_ai1_k}
\end{eqnarray}
and since the event~\(B_{k-1}\) depends only on the values of~\(\{S_j\}_{j \in {\cal C}(v_{k-1}) \cup \{v_{k-1}\}},\) we condition on~\(S_{v_k}=z\) and get from~(\ref{ai2_ai1_k}) that~\(\mathbb{P}\left(B_k(x)\right) = \int_{z=0}^{x} t_{k}(z)dF(z)\) where
\begin{eqnarray}
t_k(z) &:=& \mathbb{P}\left(\max_{j \in {\cal E}_k} S_j < z\right)  \mathbb{P}\left(B_{k-1} \bigcap \left\{ S_{v_{k-1}} < z\right\} \right) \nonumber\\
&=& \prod_{j \in {\cal E}_k} \mathbb{P}(S_j < z) \mathbb{P}\left(B_{k-1} \bigcap \left\{ S_{v_{k-1}} < z\right\} \right). \label{hawa}
\end{eqnarray}

To evaluate~(\ref{hawa}), we use the weak nested property~\(v_{k-1} \in {\cal D}(v_k) \setminus {\cal C}(v_{k-1})\) to get that~\({\cal E}_k = {\cal D}(v_k) \setminus \left({\cal C}(v_{k-1}) \cup v_{k-1}\right)\) has cardinality~\(c(v_k)-c(v_{k-1})-1.\) Consequently~
\begin{eqnarray}
t_k(z)&=& F^{c(v_k)-c(v_{k-1})-1}(z)\mathbb{P}\left(B_{k-1} \bigcap \left\{ S_{v_{k-1}} < z\right\} \right) \nonumber\\
&=& F^{c(v_k)-c(v_{k-1})-1}(z)\mathbb{P}\left(B_{k-1}(z)\right) \nonumber
\end{eqnarray}
and substituting this into the expression for~\(\mathbb{P}(B_k(x))\) determined prior to~(\ref{hawa}), we get~\[\mathbb{P}\left(B_k(x)\right) = \int_{0}^{x} F^{c(v_k)-c(v_{k-1})-1}(z)\mathbb{P}\left(B_{k-1}(z)\right) dF(z).\] By induction assumption we have~\(\mathbb{P}\left(B_{k-1}(z)\right) = \frac{F^{c(v_{k-1})}(z)}{c(v_1)\cdot c(v_2)\cdots c(v_{k-1})}\) and so
\begin{eqnarray}
\mathbb{P}\left(B_k(x)\right) &=& \frac{1}{c(v_1)\cdot c(v_2) \cdots c(v_{k-1})} \int_{0}^{x} F^{c(v_k)-1}(z)dF(z) \nonumber\\
&=& \frac{1}{c(v_1)\cdot c(v_2)\cdots c(v_k)} F^{c(v_k)}(x). \label{pbxk_up}
\end{eqnarray}
This proves the induction step and therefore completes the proof of~(\ref{homo_b}).

The first relation in~(\ref{mut_ind}) follows directly from~(\ref{homo_b}) by setting~\(x= \infty.\)  Moreover,~\(var(R_{tot}) = V_1 +2V_2\)
where
\[V_1 = \sum_{l=1}^{k} \mathbb{P}(A_{v_l}) - \left(\mathbb{P}(A_{v_l})\right)^2 = \sum_{l=1}^{k} \left(\frac{1}{d(v_l)} - \frac{1}{d^2(v_l)}\right) \leq \sum_{l=1}^{k} \frac{1}{d(v_l)} = \mathbb{E}R_{tot}\] and~\(V_2 = \sum_{1 \leq l_1 < l_2 \leq k} \left(\mathbb{P}(A_{v_{l_1}} \cap A_{v_{l_2}}) - \mathbb{P}(A_{v_{l_1}}) \mathbb{P}(A_{v_{l_2}})\right) \leq 0 \) by the first relation in~(\ref{mut_ind}). This completes the proof of~(\ref{mut_ind}).

For the mutual independence part, suppose that~\({\cal D}(v_l) \subset {\cal D}(v_{l+1})\) so that\\\(c(v_l) = d(v_l)\) for~\(1 \leq l \leq k.\) From~(\ref{homo_b}), we then get that
\begin{eqnarray}
\mathbb{P}\left(B_k(x)\right) &=& \frac{1}{d(v_1)\cdot d(v_2) \cdots d(v_{k-1})} \frac{F^{d(v_k)}(x)}{d(v_k)} \nonumber\\
&=& \left(\prod_{j=1}^{k-1} \mathbb{P}\left(A_{v_j}\right)\right)\mathbb{P}\left(A_{v_k}(x)\right) \label{mut_ind_2_dash}
\end{eqnarray}
and so the events~\(\{A_{v_l}\}_{1 \leq l \leq k-1} \cup A_{v_k}(x)\) are also mutually independent. Setting~\(x=\infty\) and using induction on~\(k,\) we then get that~\(\{A_{v_l}\}_{1 \leq l \leq k}\) are mutually independent.

Conversely suppose that~\(\{A_{v_l}\}_{1 \leq l \leq k}\) are mutually independent. By definition we know that~\(d(v_1) = c(v_1)\) and by~(\ref{mut_ind}), we have that
\[\frac{1}{c(v_1) \cdot c(v_2)} = \mathbb{P}\left(A_{v_1} \bigcap A_{v_2}\right) = \mathbb{P}(A_{v_1}) \mathbb{P}(A_{v_2}) = \frac{1}{d(v_1)} \cdot \frac{1}{d(v_2)}.\] Consequently~\(c(v_2) = d(v_2)\)
and so using
\begin{eqnarray}
c(v_2) &=& 1+\#\left({\cal D}(v_1) \cup {\cal D}(v_2)\right) \nonumber\\
&=& 1+\#\left({\cal D}(v_2)\right) + \#\left({\cal D}(v_1) \setminus {\cal D}(v_2)\right) \nonumber\\
&=& d(v_2) + \#\left({\cal D}(v_1) \setminus {\cal D}(v_2)\right),\nonumber
\end{eqnarray}
we get that~\({\cal D}({v_1}) \setminus {\cal D}(v_2) = \emptyset.\) In other words~\({\cal D}(v_1) \subseteq {\cal D}(v_2)\) and since~\(v_1 \in {\cal D}(v_2) \setminus {\cal D}(v_1),\) we get  that~\({\cal D}(v_1) \subset {\cal D}(v_2),\) strictly. The same argument repeated with the vertices~\(v_2\) and~\(v_3\) gives that~\({\cal D}(v_2) \subset {\cal D}(v_3).\) Continuing this way for~\(k\) iterations completes the proof of the Theorem.~\(\qed\)

\setcounter{equation}{0}
\renewcommand\theequation{\arabic{section}.\arabic{equation}}
\section{Constrained Colouring}\label{sec_col}
Let~\(K_n\) be the complete graph on~\(n\) vertices and let~\(X_f, f \in K_n\) be independent and identically distributed (i.i.d.)  Bernoulli random variables with
\begin{equation}\label{x_dist}
\mathbb{P}(X_f = 1)= p(f)  = 1-\mathbb{P}(X_f = 0)
\end{equation}
where~\(0 \leq p(f) \leq 1.\) Let~\(G\) be the random subgraph of~\(K_n\) formed by the set of all edges satisfying~\(X_f = 1.\) If~\(p(f)= p\) for all edges~\(f,\) then~\(G\) is said to be \emph{homogenous}; otherwise we say that~\(G\) is \emph{inhomogenous}. In this paper, we are mainly interested in colouring numbers of inhomogenous random graphs with constraints.

We begin with some well-known definitions. For an integer~\(k \geq 1,\) a proper~\(k-\)colouring of~\(G\) is a map~\(g: V \rightarrow \{1,2,\ldots,k\}\) such that~\(g(u) \neq g(v)\) for all edges~\((u,v) \in E.\) The chromatic number~\(\chi(G)\) is the smallest integer~\(k\) such that~\(G\) admits a proper~\(k-\)colouring.
It is well-known (see for example Theorem~\(2,\) Ganesan (2020)) that if~\(G\) is homogenous with edge probability~\(p(f)=p =\frac{1}{n^{\beta}}\) for some constant~\(0 < \beta < 1,\) then there are constants~\(C_1,C_2 >0\) such that
\begin{equation}\label{chi_g}
\mathbb{P}\left(\frac{C_1 n^{1-\beta}}{\log{n}} \leq \chi(G) \leq \frac{C_2n^{1-\beta}}{\log{n}}\right) = 1-o(1),
\end{equation}
where~\(o(1) \longrightarrow 0\) as~\(n \rightarrow \infty.\) 


The following result obtains upper and lower bounds for the chromatic number of an inhomogenous random graph~\(G,\) in terms of ``averaged" edge probabilities. 
\begin{Theorem}\label{thm_col} Suppose there are constants~\(0 \leq \beta_{up} \leq \beta_{low} \leq \alpha\) and positive constants~\(C_1,C_2\) such that
\begin{equation}\label{p_cond}
\frac{C_1}{n^{\beta_{low}}} \leq \frac{1}{{l \choose 2}} \sum_{u,v \in {\cal V}} p(u,v)\leq \frac{C_2}{n^{\beta_{up}}}
\end{equation}
for any set~\({\cal V} \subseteq V\) containing~\(l \geq n^{\alpha} (\log{n})^3\) vertices. There are positive constants~\(D_i,1 \leq i \leq 4\) such that
\begin{equation}\label{chi_h_bounds_low}
\mathbb{P}\left( \chi(G) \geq \frac{n^{1-\alpha}}{(\log{n})^3} \right) \geq 1-\exp\left(-D_2 \cdot n^{2\alpha-\beta_{low}} \cdot (\log{n})^6\right)
\end{equation}
and
\begin{equation}\label{chi_h_bounds}
\mathbb{P}\left(\chi(G) \leq D_3n^{\max(\alpha,1-\beta_{up})}\cdot (\log{n})^4\right) \geq 1- \exp\left(-D_4n^{2\alpha-\beta_{low}} \cdot (\log{n})^5\right)
\end{equation}
for all~\(n\) large.
\end{Theorem}

We prove the lower bound in~(\ref{chi_h_bounds}) by obtaining an upper bound for the maximum size of a stable set and prove the upper bound in~(\ref{chi_h_bounds}) using a maximum average degree estimate for the chromatic number, obtained  in Lemma~\ref{thm_comb_est} of Section~\ref{sec_comb_lem}. We provide the details below. Throughout, we use the following standard Deviation Estimate. Let~\(Z_i, 1 \leq i \leq t\) be independent Bernoulli random variables satisfying~\[\mathbb{P}(Z_i = 1) = p_i = 1-\mathbb{P}(Z_i = 0).\] If~\(W_t = \sum_{i=1}^{t} Z_i\) and~\(\mu_t = \mathbb{E}W_t,\) then for any~\(0 < \epsilon < \frac{1}{2}\) we have that
\begin{equation}\label{conc_est_f}
\mathbb{P}\left(\left|W_t-\mu_t\right| \geq \epsilon \mu_t\right) \leq 2\exp\left(-\frac{\epsilon^2}{4}\mu_t\right).
\end{equation}
For a proof of~(\ref{conc_est_f}), we refer to Corollary~\(A.1.14,\) pp.~\(312,\) Alon and Spencer (2008).

\emph{Proof of~(\ref{chi_h_bounds_low}) in Theorem~\ref{thm_col}}: For a set of vertices~\({\cal V}\) of size~\(\#{\cal V} = l,\) let
\begin{equation}\label{p_av_v}
p_{av}({\cal V}) := \frac{1}{{l \choose 2}} \sum_{u,v \in {\cal V}} p(u,v)
\end{equation}
be the average edge probability of edges formed by the vertices of~\({\cal V}.\) For a constant~\(0 < \alpha < 1\) set
\begin{equation}\label{p_low_def}
p_{low} := \min_{{\cal V}} p_{av}({\cal V}) \leq \max_{{\cal V}} p_{av}({\cal V}) =: p_{up}
\end{equation}
where the minimum and maximum are taken over all sets~\({\cal V}\) of cardinality~\(l\) satisfying~\(n^{\alpha}\cdot (\log{n})^3 \leq l \leq n.\)

If~\(\alpha(G)\) is the largest size of a stable set in the random graph~\(G,\) then we know that (Proposition~\(5.1.7,\) pp. 193, West (2000))
\begin{equation}\label{chi_alpha}
\chi(G) \geq \frac{n}{\alpha(G)}.
\end{equation}
To estimate the probability of the event~\(\{\alpha(G)\geq t\},\)  we let~\({\cal V}\) be any set of~\(t\) vertices and let~\(F({\cal V})\) be the event that~\({\cal V}\) is a stable set. We have that
\[\mathbb{P}\left(F({\cal V})\right) = \prod_{u\neq v \in {\cal V}} \left(1-p(u,v)\right) \leq e^{-{t \choose 2} \cdot p_{av}({\cal V})},\]
where~\(p_{av}({\cal V})\) is as defined in~(\ref{p_av_v}). Using~(\ref{p_cond}), we get that
\[\mathbb{P}\left(F({\cal V})\right) \leq \left(- {t \choose 2} \cdot \frac{C_1}{n^{\beta_{low}}}\right) \leq e^{-D_1t^2\cdot n^{-\beta_{low}}}\]
for some constant~\(D_1 > 0.\) Therefore letting~\(F_{tot} := \bigcup_{\cal V} F({\cal V})\) where the union is over all sets of size~\(t,\) we get by the union bound that
\begin{equation}
\mathbb{P}\left(F_{tot}\right) \leq  {n \choose t} \cdot e^{-D_1t^2\cdot n^{-\beta_{low}}} \leq  n^{t}\cdot e^{-D_1t^2\cdot n^{-\beta_{low}}}. \label{f_tot_one}
\end{equation}
For~\(t = n^{\alpha} \cdot (\log{n})^3,\) we have that
\begin{eqnarray}
n^{t}\cdot e^{-D_1t^2\cdot n^{-\beta_{low}}} &=& \exp\left(-t \log{n} \left(D_1n^{\alpha-\beta_{low}} (\log{n})^2-1\right) \right) \nonumber\\
&\leq& \exp\left(-2D_2 \cdot t n^{\alpha-\beta_{low}} \cdot (\log{n})^3\right) \nonumber\\
&=& \exp\left(-2D_2 \cdot n^{2\alpha-\beta_{low}} \cdot (\log{n})^6\right), \label{temp_q}
\end{eqnarray}
for some constant~\(D_2 > 0.\)

Substituting~(\ref{temp_q}) into~(\ref{f_tot_one}) we get that
\[\mathbb{P}(F_{tot}) \leq \exp\left(-2D_2 \cdot n^{2\alpha-\beta_{low}} \cdot (\log{n})^6\right)\] for all~\(n\) large. If the complement event~\(F^c_{tot}\) occurs, then the largest size of the stable set~\(\alpha(H) \leq n^{\alpha} \cdot (\log{n})^3.\) Plugging this into~(\ref{chi_alpha}), we get that
\[\mathbb{P}\left( \chi(H) \geq \frac{n^{1-\alpha}}{(\log{n})^3} \right) \geq 1-\exp\left(-D_2 \cdot n^{2\alpha-\beta_{low}} \cdot (\log{n})^6\right) \] and this proves~(\ref{chi_h_bounds_low}).~\(\qed\)

\emph{Proof of~(\ref{chi_h_bounds}) in Theorem~\ref{thm_col}}: We prove~(\ref{chi_h_bounds}) using the estimate~(\ref{chi_bound}) of Lemma~\ref{thm_comb_est}\((b)\) for~\(\chi(G)\) and the estimate~(\ref{max_ave_2}) for the maximum average degree~\(h_{av}(G),\) in the proof of Lemma~\ref{thm_comb_est}\((a).\)

For a set~\({\cal V} \subseteq \{1,2,\ldots,n\}\) of size~\(\#{\cal V} = l,\) let~\(m({\cal V})\) be the random number of edges in the induced subgraph~\(H_{\cal V}\) of~\(H\) formed by the vertices of~\({\cal V}.\) The expected number of edges in~\(H_{\cal V}\) is
\[\mathbb{E} m({\cal V}) = {l \choose 2}p_{av}({\cal V})\] and so by the standard deviation estimate~(\ref{conc_est_f}) we have for~\(0 < \epsilon < \frac{1}{2}\) that
\[\mathbb{P}\left(m({\cal V}) \geq {l \choose 2}p_{av}({\cal V})(1+\epsilon)\right) \leq \exp\left(-\frac{\epsilon^2}{4} {l \choose 2}p_{av}({\cal V})\right).\] We henceforth set~\(\epsilon = \frac{1}{\sqrt{\log{n}}}.\) If~\(n^{\alpha} \cdot (\log{n})^3 \leq l \leq n\) then we have from~(\ref{p_low_def}) that
\begin{equation}\label{m_av_est_one}
\mathbb{P}(m({\cal V}) \geq {l \choose 2} p_{up}(1+\epsilon)) \leq \exp\left(-\frac{\epsilon^2}{4} {l \choose 2} p_{low}\right) \leq \exp\left(-\frac{l^2\epsilon^2}{16}p_{low}\right).
\end{equation}

Defining the event
\begin{equation}
E_{bad} := \bigcup_{n^{\alpha} \leq \#{\cal V} = l \leq n} \left\{m({\cal V}) \geq {l \choose 2} p_{up}(1+\epsilon)\right\}
\end{equation}
we have from~(\ref{m_av_est_one}) and the union bound that
\begin{equation}\label{e_bad_est}
\mathbb{P}(E_{bad}) \leq \sum_{n^{\alpha} \leq l \leq n} {n \choose l} \cdot \exp\left(-\frac{l^2\epsilon^2}{16} p_{low}\right)
\end{equation}

Using~\({n \choose l}  \leq n^{l}\) we have that
\begin{equation}\label{temp_est2}
{n \choose l} \cdot \exp\left(-\frac{l^2\epsilon^2}{16} p_{low}\right) \leq n^{l} \cdot \exp\left(-\frac{l^2\epsilon^2}{16} p_{low}\right)  = \exp(-A_l),
\end{equation}
where~\(A_l = A_l(n) := l\left(\frac{l\epsilon^2\cdot p_{low}}{16} - \log{n}\right).\)
Since~\(\epsilon = \frac{1}{\sqrt{\log{n}}}\) and~\(l \geq n^{\alpha} \cdot (\log{n})^3,\) we have that~\[ l\epsilon^2\cdot p _{low} \geq n^{\alpha}  p_{low}\cdot (\log{n})^2 \geq C_1 \cdot n^{\alpha-\beta_{low}} \cdot (\log{n})^2,\] by~(\ref{p_cond}) and so using~\(\alpha \geq \beta_{low}\) and we get that
\begin{equation}\label{al_est}
A_l \geq l \cdot 2C \cdot n^{\alpha-\beta_{low}} \cdot (\log{n})^2 \geq 2C \cdot n^{2\alpha-\beta_{low}} \cdot (\log{n})^5
\end{equation}
for all~\(n\) large and some constant~\(C >0.\) Plugging~(\ref{al_est}) into~(\ref{temp_est2}) and using~(\ref{e_bad_est}) we get that
\begin{eqnarray}
\mathbb{P}(E_{bad})  &\leq& \sum_{n^{\alpha} \leq l \leq n} \exp\left(-2C n^{2\alpha-\beta_{low}}\cdot (\log{n})^5\right) \nonumber\\
&\leq& n \cdot \exp\left(-2C n^{2\alpha-\beta_{low}} \cdot (\log{n})^5\right) \nonumber\\
&\leq& \exp\left(-Cn^{2\alpha-\beta_{low}} \cdot (\log{n})^5\right) \label{e_bad_est2}
\end{eqnarray}
for all~\(n\) large. 

If the complement event~\(E^c_{bad}\) occurs, then in any induced subgraph~\(H_l\) of~\(H\) containing~\(T := n^{\alpha} \cdot (\log{n})^3 \leq l \leq n\) vertices, the number of edges is at most~\(m(H_l) \leq {l \choose 2} p_{up}(1+\epsilon)\) and so
\begin{equation}\label{cruc}
\max_{T+1 \leq l \leq n} \frac{2m(H_l)}{l} \leq (n-1) \cdot  p_{up} (1+\epsilon) \leq C_2 \cdot n^{1-\beta_{up}}(1+\epsilon)
\end{equation}
where~\(C_2 > 0\) is the constant in~(\ref{p_cond}). Substituting~(\ref{cruc}) into~(\ref{max_ave_2}) and using the fact that~\(\epsilon = \frac{1}{\sqrt{\log{n}}},\) we get that the maximum average degree
\begin{equation}
h_{av} \leq n^{\alpha} \cdot (\log{n})^3 + C_2 n^{1-\beta_{up}}(1+\epsilon) \leq D \cdot n^{\max(\alpha,1-\beta_{up})} \cdot (\log{n})^3 \label{h_av_est}
\end{equation}
with probability at least~\(1-\omega_n,\) for some constant~\(D > 0.\) Here~\(\omega_n\) is the final expression in~(\ref{e_bad_est2}). Plugging~(\ref{h_av_est}) into~(\ref{chi_bound}) of Lemma~\ref{thm_comb_est}\((b),\) we get the upper bound in~(\ref{chi_h_bounds}). This completes the proof of Theorem~\ref{thm_col}.~\(\qed\)

\subsection*{Weighted Colouring}
In this subsection, we study a weighted generalization of proper colouring defined as follows. Equip each edge~\(h\) of the random graph~\(G\) obtained in~(\ref{x_dist}) with a deterministic weight~\(w(h) \geq 1.\) Let~\(f : V \rightarrow \{1,2,\ldots\}\) be any map from the vertex set of~\(H\) to the set of positive integers.
\begin{Definition} We say that~\(f\) is a proper~\(w-\)\emph{weighted colouring} or simply proper weighted colouring of~\(H\) if for each edge~\(h = (u,v) \in E\) with endvertices~\(u\) and~\(v\) we have
\begin{equation}\label{rad_cond}
|f(u)-f(v)| \geq w(u,v).
\end{equation}
We define the weighted colouring number of~\(H\) as
\begin{equation}\label{grl_def}
\chi_{w} = \chi_{w}(H)  := \min_{f} \max_{(u,v) \in E} |f(u)-f(v)|,
\end{equation} where the minimum is over all proper weighted colourings of~\(H.\) 
\end{Definition}
We could interpret~\(\chi_{w}(H)\) as a measure of the effect of weights on a proper colouring of~\(H.\) Indeed if~\(w(u,v) =1\) for all edges, then the weighted colouring number~\(\chi_{w}(H) = \chi(H),\) the chromatic number. In general,~\(\chi_{w}(H) \geq \chi(H)\) and if~\(w(u,v) \leq K\) for some~\(K\) and all edges~\((u,v),\) then~\(\chi_{w}(H) \leq K \chi(H).\)

Henceforth we consider weighted colouring number of~\(G\) when the edge weights are random and unbounded. Formally, we equip the~\({n \choose 2}\) edges of~\(K_n\) with independent and identically distributed (i.i.d.) weights~\(\{w(h)\}_{h \in K_n}\) having a constant mean~\(1 \leq \mu_0 < \infty\) and satisfying~\(w(h) \geq 1\) a.s.\ for all edges~\(h \in K_n.\) The edge weights are also independent of~\(G.\)



For a realization~\(\omega\) of the random graph~\(G,\) let~\(\mathbb{P}_{\omega}\) be the probability measure associated with the edge weights of~\(\omega.\) For constants\\\(C,\epsilon_1,\epsilon_2,\epsilon_3 >0\) we define the event
\[E_{good}(C,\epsilon_1,\epsilon_2) := \left\{\omega : \mathbb{P}_{\omega}\left(\chi_w(\omega) \leq C n^{\epsilon_1}\right) \geq 1-\frac{C}{n^{\epsilon_2}} \right\} \] and have the following result regarding the weighted colouring number~\(\chi_{w}(G).\) For a vertex~\(u,\) we define of the average edge probability~\(p_{av}(u) := \frac{1}{n-1}\sum_{v \neq u} p(u,v).\)
\begin{Theorem}\label{thm_grl} Suppose~\(\mathbb{E}w^{2s}(u,v) <\infty\) for some integer~\(s \geq 1.\) \\
\((a)\) Suppose for each vertex~\(u\) the average edge probability satisfies~\(\frac{C_1}{n^{\beta}} \leq p_{av}(u) \leq \frac{C_2}{n^{\beta}}\) for some positive constants~\(C_1,C_2\) and~\(0 < \beta < 1.\) For every~\(s > \frac{1}{1-\beta}\) there is a constant~\(C > 0\) such that
\begin{equation}\label{rgg_chi_w}
\mathbb{P}\left(E_{good}(C,1-\beta,s(1-\beta)-1)\right) \geq 1-e^{-C^{-1}n^{1-\beta}}.
\end{equation}
\((b)\) Suppose the average edge probability~\(p_{av} := {n \choose 2}^{-1} \sum_{h} p(h)\) satisfies~\(\frac{C_1}{n^{\beta}} \leq p_{av} \leq \frac{C_2}{n^{\beta}}\) for some~\( \frac{2}{s} < \beta  < 1.\)   For every~\(0 < \gamma < \frac{\beta}{2}-\frac{1}{s},\) there is a constant~\(D>0\) such that
\begin{equation}\label{e_good}
\mathbb{P}(E_{good}(D,1-\gamma,\beta)) \geq 1- e^{-D^{-1}n^{2-\beta}}
\end{equation}
\end{Theorem}
Part~\((a)\) of the above result essentially says that for nearly all realizations, the weighted colouring number is at most of the order of~\(n^{1-\beta}\) with high probability. From part~\((b)\) of Theorem~\ref{thm_grl}, we see that even under the weaker condition on the overall edge probability average, the weighted colouring number is still~\(o(n)\) with high probability, provided the weights have a bounded~\(s^{th}\) moment for~\(s\) sufficiently large. Here and henceforth~\(o(n)\) denotes a sequence satisfying~\(\frac{o(n)}{n} \longrightarrow 0\) as~\(n \rightarrow \infty.\)






\emph{Proof of Theorem~\ref{thm_grl}\((a)\)}: For~\(1 \leq i \leq n\) let~\(N_i\) be the number of neighbours of the vertex~\(i\) in the graph~\(G.\) From the bounds on the average edge probability~\(p_{av}(i),\) we see that the expected number of neighbours for vertex~\(i\) lies~\(C_1n^{1-\beta}\) and~\(C_2n^{1-\beta}.\) Therefore by the concentration estimate~(\ref{conc_est_f}) with~\(\epsilon = \frac{1}{2},\) we have that
\[\mathbb{P}\left(\frac{C_1n^{1-\beta}}{2} \leq N_i \leq 2C_2 n^{1-\beta}\right) \geq 1-e^{-Dn^{1-\beta}},\] for some positive constant~\(D.\) Letting~\(E_{nei} :=\bigcap_{1 \leq i \leq n} \{\frac{C_1 n^{1-\beta}}{2} \leq N_i \leq 2C_2 n^{1-\beta}\},\) we have by the union bound that
\begin{equation}\label{nei_est}
\mathbb{P}(E_{nei}) \geq 1- ne^{-D n^{1-\beta}}.
\end{equation}

Let~\(\omega \in E_{nei}\) be a realization of~\(G.\) To find an upper bound for~\(\chi_{w}(\omega),\) we use the locally averaged bound~(\ref{rk_est_det2}) in Proposition~\ref{grl_lem}. For~\(1 \leq i \leq n\) let~\(J_i := \sum_{v \sim i} w(i,v)\) be the sum of the weights of edges containing~\(i\) as an endvertex. Using the fact that there are at most~\(2C_2n^{1-\beta}\) nodes adjacent to~\(i,\) we have that~\(J_i\) is stochastically dominated by~\(\sum_{i=1}^{2C_2n^{1-\beta}} Z_i\) where~\(Z_i\) are i.i.d. with the same distribution as the edge weights. Using the fact that the edge weights have bounded~\(2s^{th}\) moment, we now show that there are constants~\(D_1,D_2 > 0\) such that
\begin{equation}\label{legally_blonde}
\mathbb{P}\left(\sum_{i=1}^{2C_2 n^{1-\beta}} Z_i \geq D_1 n^{1-\beta}\right) \leq \frac{D_2}{n^{s(1-\beta)}}.
\end{equation}
Indeed, denoting the expectation of~\(Z_i\) as~\(\mu_Z\) and using the Markov inequality we have for any integer~\(t \geq 2s\) that
\begin{equation}\label{leg_blonde2}
\mathbb{P}\left(\sum_{i=1}^{t}Z_i \geq 2t \mu_Z\right) \leq \frac{\mathbb{E}\left(\sum_{i=1}^{t}(Z_i-\mu_Z)\right)^{2s}}{(t\mu_Z)^{2s}}.
\end{equation}
Expanding~\(\left(\sum_{i=1}^{t}(Z_i-\mu_Z)\right)^{2s},\) we see that any term of the form\\\(\prod_{j}(Z_{i_j}-\mu_Z)^{b_j}\) has expectation zero, unless each~\(b_j\) is even. This implies that the total number of terms with non-zero expectation is~\(\sum_{l=1}^{s}{t \choose l} \leq s {t \choose s} \leq st^{s}\) by the monotonicity of the binomial coefficient for~\(t \geq 2s.\) Thus
\[\mathbb{E}\left(\sum_{i=1}^{t}(Z_i-\mu_Z)\right)^{2s} \leq C_1 s t^{s}\] for some constant~\(C_1 > 0\) and plugging this into~(\ref{leg_blonde2}), we get that\\\(\mathbb{P}\left(\sum_{i=1}^{t}Z_i \geq 2t \mu_Z\right) \leq \frac{C_2}{t^{s}}.\) Setting~\(t = 2C_2n^{1-\beta}\) then gives us~(\ref{legally_blonde}).

Letting~\(F_{wt} := \bigcap_{1 \leq i \leq n}\{J_i \leq D_1n^{1-\beta}\}\) and using the union bound, we get from~(\ref{legally_blonde}) that
\begin{equation}\label{legally_blonde3}
\mathbb{P}_{\omega}\left(F_{wt}\right) \geq 1-\frac{D_2}{n^{s(1-\beta)-1}}.
\end{equation}
If~\(F_{wt}\) occurs, then using~(\ref{rk_est_det2}), we get that~\(\chi_{w} \leq D_3 n^{1-\beta}\) for some constant~\(D_3 > 0.\) Thus we get~(\ref{rgg_chi_w}) from~(\ref{nei_est}) and~(\ref{legally_blonde3}).~\(\qed\)


\emph{Proof of Theorem~\ref{thm_grl}}\((b)\): The expected number of edges in~\(G\) equals~\({n \choose 2} p_{av}\) and so if~\(E_{edge}\) is the event that the number of edges in the random graph~\(G\) lies between~\({n \choose 2}\frac{p_{av}}{2}\) and~\({n \choose 2}\frac{3p_{av}}{2},\) then by the concentration estimate~(\ref{conc_est_f}) with~\(\epsilon = \frac{1}{2},\) we have that
\begin{equation}\label{e_edge_est}
\mathbb{P}(E^c_{edge}) \leq \exp\left(-C_1{n \choose 2}p_{av}\right) \leq \exp\left(-C_2n^{2-\beta}\right)
\end{equation}
for some constants~\(C_1,C_2 >0.\)

Let~\(\omega \in E_{edge}\) be any realization of~\(G\) and let~\(q = q(\omega)\) be the number of edges in~\(\omega.\) By definition~\(q\) is at least of the order of~\(n^{2-\beta}\) and so  if~\(W_{tot}\) denotes the total weight of edges in~\(\omega,\) then by Chebychev's inequality, we have for constant~\(\epsilon > 0\) that
\begin{equation}\label{w_tot_est}
\mathbb{P}_{\omega}(|W_{tot}-q\mu_0| \geq \epsilon q\mu_0) \leq \frac{var(W_{tot})}{\epsilon^2q^2\mu_0^2} \leq \frac{C_3}{q} \leq \frac{C_4}{n^{2-\beta}}
\end{equation}
for some constants~\(C_3,C_4>0.\) Next let~\(E_{wt}\) be the event that the maximum edge weight in~\(\omega\) is at most~\(n^{(2+\delta)/s}\) for some constant~\(\delta > 0\) to be determined later. By the bounded~\(s^{th}\) moment assumption of the edge weights and the Markov inequality we have for any edge~\(h\) of~\(\omega\) that
\[\mathbb{P}_{\omega}\left(w(h) \geq n^{(2+\delta)/s}\right) \leq \frac{C_5}{n^{2+\delta}}\] for some constant~\(C_5 > 0.\) Since the number of edges of~\(\omega\) is at most of the order of~\(n^{2-\beta},\) we get by the union bound that
\begin{equation}\label{e_wt_est}
\mathbb{P}_{\omega}(E^c_{wt}) \leq \frac{C_6}{n^{\beta+\delta}}
\end{equation}
for some constant~\(C_6 > 0.\)

Set~\(\epsilon = 1\) and suppose that~\(E_{tot} := E_{wt} \cap \{W_{tot} \leq 2q\mu_0\}\) occurs. From~(\ref{w_tot_est}) and~(\ref{e_wt_est}), we get that
\begin{equation}\label{e_tot_est}
\mathbb{P}_{\omega}(E_{tot}) \geq 1- \frac{C_7}{n^{2-\beta}} - \frac{C_7}{n^{\beta+\delta}}
\end{equation}
for some constant~\(C_7 > 0.\) We choose~\(\delta > 0\) small so that~\(2-\beta >\beta+\delta\) and this is possible since~\(\beta < 1.\)  If the event~\(E_{tot}\) occurs, then setting~\(K = n^{(2+\delta)/s}\) and~\(m\mu = W_{tot}\)  in the bound~(\ref{rk_est_det}), we get that~\[\chi_{w}(\omega) \leq C_8 n^{\frac{2+\delta}{2s}} \cdot \sqrt{q} \leq C_9 n^{\frac{2+\delta}{2s}} n^{1-\frac{\beta}{2}}= C_9 n^{1-\left(\frac{\beta}{2}-\frac{1}{s} - \frac{\delta}{2s}\right)}\]
for some constants~\(C_8,C_9 > 0.\) From the statement of the theorem we know that~\(\beta  > \frac{2}{s}\) and so we choose~\(\delta > 0\) smaller if necessary so that~\(\frac{\beta}{2}-\frac{1}{s} - \frac{\delta}{2s}  >0\) as well. With this choice of~\(\delta\) and setting~\(E_{good} = E_{edge},\) we obtain~(\ref{e_good}) from~(\ref{e_edge_est}) and~(\ref{e_tot_est}).~\(\qed\)

\setcounter{equation}{0}
\renewcommand\theequation{\arabic{section}.\arabic{equation}}
\section{Combinatorial Lemmas}\label{sec_comb_lem}
In this section, we state and prove the combinatorial lemmas used in the proof of the Theorems in the previous section. Throughout~\(H =(V,E)\) is a deterministic graph with vertex set~\(V = \{1,2,\ldots,n\}\) and edge set~\(E.\)

Our first result in concerns an upper bound for the chromatic number of a graph in terms of its maximum average degree, used in the proof of Theorem~\ref{thm_col}. We begin with a couple of definitions. For a graph~\(H =(V,E),\) we define the maximum average degree as~\(h_{av} = h_{av}(H) := \max_{\Gamma} d_{av}(\Gamma),\) where the maximum is over all subgraphs~\(\Gamma \subseteq G\) and~\(d_{av}(\Gamma)\) is the average vertex degree in~\(\Gamma.\) By definition, we have that~\(d_{av} \leq h_{av} \leq \Delta,\) the maximum degree of a vertex in~\(H.\) The following result obtains bounds for the chromatic number of~\(H\) in terms of its maximum average degree.
\begin{Lemma}\label{thm_comb_est} For any graph~\(H\) containing~\(n\) vertices and~\(m\) edges, we have that the chromatic number
\begin{equation}\label{chi_bound}
\chi(H) \leq 2\max(1,h_{av})\cdot \log\left(\frac{ne}{\max(1,h_{av})}\right),
\end{equation}
where~\(h_{av}\) satisfies~\(d_{av} \leq h_{av} \leq \min\left(\sqrt{2m},6(m\Delta)^{1/3}\right).\) 
\end{Lemma}
From the above Lemma, we see that if~\(H\) has bounded maximum average degree, then~\(\chi(H) = O(\log{n}).\)  For context we recall from the edge count bound that if the average degree of~\(H\) is bounded, then~\(\chi(H) \leq 2\sqrt{m}= \sqrt{2d_{av}(G)} \cdot \sqrt{n} = O(\sqrt{n}).\)   This is the best possible since if~\(H\) contains a complete subgraph formed by~\(\sqrt{n}\) vertices, then the chromatic number is at least~\(\sqrt{n}.\) On the other hand, if~\(H\) has bounded maximum average degree, then \emph{every} subgraph of~\(H\) is also sparse and so we expect~\(H\) to have low chromatic number. This is reflected in the estimate~(\ref{chi_bound}).


\emph{Proof of Lemma~\ref{thm_comb_est}}:  A stable set in~\(H\) is a set of vertices no two of which are adjacent in~\(H.\) Let~\({\cal I}_1 = \{v_1,\ldots,v_t\}\) be a maximum stable set, i.e. a stable set of maximum size, in~\(G_1 := H.\) We assign the colour~\(1\) to all vertices in~\({\cal I}_1.\) We remove all the vertices of~\({\cal I}_1\) from~\(G_1\) to obtain a graph~\(G_2.\) We now repeat the above procedure with~\(G_2.\) Letting~\({\cal I}_2\) be the maximum stable set in~\(G_2,\) we assign the colour~\(2\) to all vertices in~\({\cal I}_2.\) Continuing the above procedure, let~\(G_{k+1}\) be the graph obtained at the end of~\(k\) iterations. The partial colouring of~\(H\) obtained so far uses~\(k\) colours and by construction is proper. Therefore if~\(G_{k+1}\) has~\(n_{k+1}\) vertices, then~\(\chi(G_{k+1}) \leq n_{k+1}\) and so
\begin{equation}\label{chi_g_k}
\chi(G) \leq k + n_{k+1}.
\end{equation}

To estimate~\(n_{k+1},\) we estimate the size of~\({\cal I}_i\) for each~\(1 \leq i  \leq k.\) Recalling that~\(d_{av}\) and~\(h_{av} \geq d_{av}\) are the average degree and maximum average degree of~\(H,\) respectively, we get from Theorem~\(3.2.1,\) pp.~\(29,\) Alon  and Spencer (2008), that~\({\cal I}_1\)  has cardinality
\begin{equation}\label{alp_def}
\#{\cal I}_1 \geq  \frac{n}{2\max(1,d_{av})} \geq \frac{n}{2\max(1,h_{av})} := n\alpha.
\end{equation}
From~(\ref{alp_def}), we see that the graph~\(G_2\) has~\(n_2 \leq n(1-\alpha)\) vertices. The graph~\(G_2\) has a maximum average degree of at most~\(h_{av}\) and so arguing as in~(\ref{alp_def}) we get that~\(\#{\cal I}_2 \geq n_2 \alpha.\) This implies that~\(G_3\) has~\(n_3 \leq n_2(1-\alpha) \leq n(1-\alpha)^2\) vertices.  Continuing this way, we get that~\(G_{k+1}\) has~\(n_{k+1} \leq n(1-\alpha)^{k}\) vertices and  so from~(\ref{chi_g_k}), we get that~\(\chi(G) \leq k + ne^{-\theta k} =:g(k)\) where~\(\theta := |\log(1-\alpha)|.\)

We see that~\(g(k)\) is minimized if~\(k\) satisfies~\(1- n\theta \cdot e^{-\theta k} = 0\) and for this value of~\(k,\) we get that~\(\chi(G) \leq \frac{\log(n\theta)+1}{\theta}.\) Using~\(|\log(1-x)| > x\) we get that~\(\theta > \alpha\) and since~\(\alpha = \frac{1}{2\max(1,h_{av})} \leq \frac{1}{2},\) we also get that~\[\theta  = | \log(1-\alpha)| \leq \sum_{j \geq 1} \alpha^{j} = \frac{\alpha}{1-\alpha} \leq 2\alpha.\] In effect~\(\alpha \leq \theta \leq 2\alpha\)
and so \[\chi(G) \leq \frac{\log(2\alpha\theta)+1}{\alpha} = 2\max(1,h_{av})\cdot \left(\log\left(\frac{n}{\max(1,h_{av})}\right) + 1\right).\] This proves~(\ref{chi_bound}).

To obtain the bounds for~\(h_{av},\) we use the fact that for any integer~\(1 \leq T \leq n,\) the maximum average degree satisfies
\begin{eqnarray}\label{max_ave_2}
h_{av} &=& \max\left(\max_{1 \leq l \leq T} \max_{H_l} d_{av}(H_l), \max_{T+1 \leq l \leq n} \max_{H_l} d_{av}(H_l)\right) \nonumber\\
&\leq& \max\left(T, \max_{T+1 \leq l \leq n} \max_{H_l} d_{av}(H_l)\right) \nonumber\\
&=& \max\left(T,\max_{T+1 \leq l \leq n} \max_{H_l} \frac{2m(H_l)}{l}\right)
\end{eqnarray}
where~\(H_l\) is an \emph{induced} subgraph of~\(G\) containing~\(l\) vertices and~\(m(H_l)\) is the number of edges in~\(H_l.\)  For any subgraph~\(H_l\) we have that~\(m(H_l) \leq m,\) the total number of edges in~\(H\) and so from~(\ref{max_ave_2}) we see that
\begin{equation}\label{h_avf}
h_{av}(G) \leq \max\left(T, \max_{T+1 \leq l \leq n} \frac{2m}{l}\right) \leq \max\left(T,\frac{2m}{T}\right)
\end{equation}
and  the final expression in~(\ref{h_avf}) attains its maximum at~\(T= \sqrt{2m}.\) Thus we get that~\(h_{av}(G) \leq \sqrt{2m}.\)

For the final estimate~\(h_{av} \leq (m\Delta)^{1/3},\) we again use~(\ref{max_ave_2}) as follows. For~\(x > 0,\) the probability that a randomly chosen vertex has degree larger than~\(x\) is bounded above by the Markov inequality as~\(\frac{d_{av}}{x}\) and so the number of vertices in~\(H\) with degree larger than~\(x\) is at most~\(\frac{nd_{av}}{x} = \frac{m}{x}.\) For any subgraph~\(H_l\) on~\(l\) vertices, we therefore have that the number of edges~\(m(H_l) \leq \frac{m}{x} \cdot \Delta + l x\) and plugging this into~(\ref{max_ave_2}), we get that
\begin{eqnarray}
h_{av} &\leq& \max\left(T, \max_{T+1 \leq l \leq n} \frac{2m\Delta}{lx} + 2x\right) \nonumber\\
&\leq& \max\left(T, \frac{2m\Delta}{Tx}+ 2x\right) \nonumber\\
&\leq& T + \frac{2m\Delta}{Tx} + 2x. \label{poojax}
\end{eqnarray}
The final expression in~(\ref{poojax}) is minimized at~\(T = \sqrt{\frac{2m\Delta}{x}}\) and so
\[h_{av} \leq 2\sqrt{\frac{2m\Delta}{x}} + 2x.\] Setting~\(x = (m\Delta)^{1/3},\) we then get the desired bound~\(h_{av} \leq (2\sqrt{2}+2)(m\Delta)^{1/3}\) and this completes the proof of the Lemma.~\(\qed\)

To prove Theorem~\ref{thm_grl}, we use the following deterministic result that obtains bounds for the weighted colouring number~\(\chi_{w}(H)\) in terms of averaged and locally averaged edge weights.
\begin{Lemma} \label{grl_lem} Suppose~\(H\) is a graph with~\(m\) edges and let~\(w(u,v)\) be the weight associated with edge~\((u,v).\)\\
\((a)\) We have that
\begin{equation}\label{rk_est_det2}
\chi_{w}(H) \leq 1+\max_{v \in H} \sum_{u \sim v} (2w(u,v)-1).
\end{equation}
\((b)\) If~\(w(h) \leq K\) for some~\(K > 0\) and all edges~\(h \in H,\) then
\begin{equation}\label{rk_est_det}
\chi_{w}(H) \leq 2\sqrt{2mK\mu}
\end{equation}
where~\(\mu := \frac{1}{m}\sum_{h}w(h)\) is the average edge weight in~\(H.\)
\end{Lemma}
The bound in part~\((a)\) of Lemma~\ref{grl_lem} is a generalization of the maximum degree bound~\(\chi(H) \leq \Delta+1.\) In fact if~\(w(u,v) = 1\) for all edges~\((u,v) \in H,\) then the right side of~(\ref{rk_est_det2}) reduces to~\(\Delta+1.\) As in the proof of the maximum degree bound, we use a greedy colouring procedure to obtain~(\ref{rk_est_det2}). Part~\((b)\) of Proposition~\ref{grl_lem} establishes bounds for the weighted colouring number in terms of the average edge weight. The bound in~(\ref{rk_est_det}) is essentially the best possible since if~\(H\) contains a clique of size~\(c \sqrt{m}\) for some~\(c > 0\) then using the fact that each edge weight is at least one, we get~\(\chi_{w}(H) \geq c\sqrt{m}.\)

We also remark here that a weaker upper bound of~\(\chi_w(H) \leq 2K\sqrt{m}\) can be obtained using a direct counting argument as follows: Let~\({\cal C} := \{K,2K,\ldots,rK\},\) where~\(r = \chi(H)\) is the smallest possible integer that allows a proper colouring of~\(G\) using only  colours from~\({\cal C}.\) Let~\({\cal V}_i\) be the set of all vertices with colour~\(iK\) for~\(1 \leq i \leq r.\) By the minimality of~\(r,\) there must exist at least one edge between a vertex in~\({\cal V}_i\) and a vertex in~\({\cal V}_j,\) for any~\(i \neq j\) and so~\({r \choose 2} \leq m.\) Consequently~\(\chi_w \leq rK \leq 2K \sqrt{m}.\) For~\(\mu \leq \frac{K}{3}\) the bound in~(\ref{rk_est_det}) is stronger and is obtained using the probabilistic method below.

We now prove Lemma~\ref{grl_lem} and Theorem~\ref{thm_grl} in that order below.\\
\emph{Proof of Lemma~\ref{grl_lem}\((a)\)}: Let~\(L +1:= 1+\max_{v \in H} \sum_{u \sim v} (2w(u,v)-1),\) the right hand side of~(\ref{rk_est_det2}). We iteratively colour the vertices of~\(H\) from the set~\(\{1,2,\ldots,L+1\}.\) Pick a vertex~\(u_1\) and assign the colour~\(l(u_1) = 1.\) For the~\(i^{th}\) iteration,~\(i \geq 2,\) we pick an uncoloured vertex~\(u_i\) and let~\(l(v_1),\ldots,l(v_t)\) be the colours of its coloured neighbours. Let~\(I(v_j)\) be the set of all integers~\(l\) satisfying~\(|l-l(v_j)| \leq w(u_i,v_j)-1.\) There are~\(2w(u_i,v_j)-1\) integers in~\(I(v_j)\) and so the number of integers in~\(\bigcup_{1 \leq j \leq t} I(v_j)\) is at most~\(L.\) We therefore assign a colour~\(l(u_i) \in \{1,2,\ldots,L+1\} \setminus \bigcup_{1 \leq j \leq t}I(v_j)\) to the vertex~\(u_i.\) Continuing this process we obtain a proper weighted colouring of~\(H\) and this proves Lemma~\ref{grl_lem}\((a).\)~\(\qed\)

\emph{Proof of Lemma~\ref{grl_lem}\((b)\)}: For integer~\(\theta \geq 1\) to be determined later, let~\(X_i, 1 \leq i \leq n\) be independent and identically distributed (i.i.d.) random variables uniformly randomly chosen from~\(\{1,2,\ldots, \theta\}.\) Say that edge~\(h = (u,v)\) is bad if~(\ref{rad_cond}) does not hold. Let~\(h_i = (u_i,v_i), 1 \leq i \leq q\) be the set of all bad edges. We now iteratively change the labels assigned to the vertices of bad edges until no bad edge remains. We begin by setting~\(w_1 := u_1\) and~\(f(w_1) = \theta+K\) and marking~\(u_1.\) We then pick an unmarked endvertex~\(w_2\) of some bad edge and set~\(f(w_2) = \theta+2K.\) We proceed this way until each endvertex of a bad edge is marked. Finally, for the rest of the vertices we set~\(f(u) = X_u.\)

Since~\(w(h) \leq K\) for all~\(h,\) we have that~\(f\) is a proper weighted colouring of~\(H.\) Moreover, if the number of bad edges is~\(N_{bad},\) then~\(1 \leq f(u) \leq \theta + N_{bad}K\) and so
\begin{equation}\label{rk_est_ax}
\chi_{w} \leq \theta + N_{bad} K.
\end{equation}
To bound~\(N_{bad},\) we first estimate~\(\mathbb{E}N_{bad}\) as follows. For any two vertices~\(u\) and~\(v\) and any integer~\(1 \leq j \leq \theta-1,
\) we have that~\(|X_u-X_v| = j\) if and only if either~\(X_u = l\) and~\(X_v = l+j\) for some~\(1 \leq l \leq \theta-j\) or~\(X_v = l\) and~\(X_u = l+j\) for some~\(1 \leq l \leq \theta-j.\) Thus~\(\mathbb{P}(|X_u-X_v| = j) = \frac{2}{\theta^2}(\theta-j)\) and so for any~\(t \leq \theta,\) we have that
\[\mathbb{P}(|X_u-X_v| \leq t) = \frac{2}{\theta^2}\left(\theta t- \frac{t(t+1)}{2}\right) \leq \frac{2t}{\theta}.\] For an edge~\(h\) with weight~\(w(h),\) we have that the probability that~\(h\) is bad is at most~\(\frac{2w(h)}{\theta}\) and so
\begin{equation}
\mathbb{E}N_{bad} \leq \sum_{h} \frac{2w(h)}{\theta} = \frac{2m\mu}{\theta} \label{en_bad}
\end{equation}
where~\(\mu = \frac{1}{m}\sum_{h}w(h)\) is the average edge weight.

From~(\ref{en_bad}), we see that there exists a realization such that~\(N_{bad} \leq \frac{2m\mu}{\theta}\) and plugging this into~(\ref{rk_est_ax}), we get that
\begin{equation}\label{rk_est2}
\chi_{w} \leq \theta + \frac{2mK\mu}{\theta}
\end{equation}
and the right side of~(\ref{rk_est2}) is minimized at~\(\theta = \sqrt{2mK\mu}.\) Thus~\(\chi_{w}(H) \leq 2\sqrt{2mK\mu}\) and this completes the proof of the Lemma.~\(\qed\)

\subsection*{Acknowledgement}
I thank Professors Rahul Roy, Alberto Gandolfi, C. R. Subramanian and K. Adhikari for crucial comments and also thank IMSc and IISER Bhopal for my fellowships.

\bibliographystyle{plain}

\end{document}